\documentclass[12pt]{amsart}
\usepackage{amsmath,amsfonts,amssymb,amsrefs,comment,url,amscd}
\usepackage[usenames]{color}
\usepackage{a4wide}
\usepackage{hyperref}

\newcommand{\Z}{\mathbb{Z}}
\newcommand{\C}{\mathbb{C}}
\newcommand{\m}{\mathfrak{m}}
\newcommand{\n}{\mathfrak{n}}
\newcommand{\abs}[1]{\left|{#1}\right|}
\newcommand{\diag}{\operatorname{diag}}
\newcommand{\GL}{\operatorname{GL}}
\newcommand{\Aut}{\operatorname{Aut}}

\begin{document}

\title[Explicit decomposition of certain induced representations]
{Explicit decomposition of certain induced representations of the general linear group}
\author{Erez Lapid}
\date{}
\address{Department of Mathematics, Weizmann Institute of Science, Rehovot 7610001, Israel}
\email{erez.m.lapid@gmail.com}
\maketitle

\begin{abstract}
We provide some experimental results on the decomposition of the parabolic induction of $\pi\otimes\pi$
in the Grothendieck group where $\pi$ is an irreducible representation of $\GL_n$ over a local non-archimedean field.
\end{abstract}

The purpose of this short note is to report on some experimental results in the representation theory
of the general linear group over a non-archimedean local field $F$.
This theory was developed in the fundamental work of Bernstein--Zelevinsky \cites{MR0425031, MR0579172}.

The direct sum over $n\ge0$ of the categories of finite length, smooth, complex representations of $\GL_n(F)$,
forms a ring category\footnote{i.e., a locally finite abelian monoidal category with bilinear and biexact tensor product}
$\mathcal{C}$ with respect to parabolic induction, which we denote by $\times$.
(The unit object is the one-dimensional representation of the trivial group $\GL_0(F)$.)

For any $a,b\in\Z$ with $a\le b$ we call $[a,b]$ a \emph{segment} and denote by $Z([a,b])$ the character
$\abs{\det}^{(a+b)/2}$ of $\GL_{b-a+1}(F)$.
By definition, a \emph{multisegment} is a formal finite sum of segments $\m=\sum_{i=1}^k\Delta_i$.
The set of multisegments forms a commutative monoid.
Writing $\Delta_i=[a_i,b_i]$ with $a_1\ge\dots\ge a_k$, the induced representation
\[
Z(\Delta_1)\times\dots\times Z(\Delta_k)
\]
admits an irreducible socle which depends only on $\m$, up to equivalence.
We denote it by $Z(\m)$.
By Zelevinsky classification \cite{MR584084}, the map $\m\mapsto Z(\m)$ is injective. Its image consists
of the irreducible representations of $\GL_n(F)$, $n\ge0$ whose Jacquet module with respect to the Borel subgroup
contains a character of the diagonal torus of the form
$\diag(t_1,\dots,t_n)\mapsto\prod_{i=1}^n\abs{t_i}^{\lambda_i}$ where $\lambda_i\in\Z$ for all $i$.
One of the nice properties of this classification is that for any multisegments $\m_1$ and $\m_2$,
the representation $Z(\m_1+\m_2)$ occurs with multiplicity one in the Jordan--H\"older sequence
of $Z(\m_1)\times Z(\m_2)$. In particular, if $Z(\m_1)\times Z(\m_2)$ is irreducible, then it is necessarily equivalent
to $Z(\m_1+\m_2)$.

The irreducible representations of the form $Z(\m)$ generate a Serre subcategory\footnote{i.e., a full subcategory
closed under subobjects, quotients and extensions} $\mathcal{C'}$ of $\mathcal{C}$
that is closed under parabolic induction.
In fact, one can decompose $\mathcal{C}$ as a direct sum of categories, each of which is
equivalent (as a ring category) to $\mathcal{C'}$.
Thus, as far as understanding the decomposition of parabolic induction of objects in $\mathcal{C}$ is concerned,
the restriction to $\mathcal{C'}$ entails no loss in generality.


Let $V$ be a finite-dimensional $\Z$-graded vector space over $\C$.
The group $\Aut(V)$ of grading preserving automorphisms of $V$ acts with finitely many orbits on the space
$E_+(V)$ of degree $1$ (nilpotent) graded endomorphisms of $V$.
By a graded version of the Jordan normal form, the orbits are indexed by multisegments whose support (in the obvious sense)
is the graded dimension of $V$ \cite{MR617466}.
Similarly, we may consider the space $E_-(V)$ of degree $-1$ graded endomorphisms of $V$.
It is the dual space to $E_+(V)$ under the trace pairing, which is $\Aut(V)$-invariant.
By a general principle \cite{MR0390138}, the orbits in $E_+(V)$ (as well as those of $E_-(V)$)
correspond bijectively to the irreducible components of the commuting variety
\[
\mathcal{X}(V)=\{(A,B)\in E_+(V)\times E_-(V):AB=BA\}.
\]
We denote by $\mathcal{C}_{\m}$ the irreducible component corresponding to (the orbit indexed by) a multisegment $\m$.
Note that $\Aut(V)$ acts on $\mathcal{X}(V)$ and preserves its irreducible components.

An interesting class of irreducible representations of $\GL_n(F)$ are those $\pi$ for which $\pi\times\pi$ is irreducible.
They were called \emph{$\square$-irreducible} in \cite{MR3866895} and studied there following the work
of Kang--Kashiwara--Kim--Oh \cites{MR3314831, MR3758148} (in a broader context which will not be discussed here).
In particular, the reducibility properties of $\pi_1\times\pi_2$ are more tractable if at least one of $\pi_1$ or $\pi_2$
is $\square$-irreducible. (We refer to \cite{MR3866895} for more details.)

In their work (again, in a much more general context), Gei{\ss}--Leclerc--Schr{\"o}er  highlighted the property
that $\mathcal{C}_{\m}$ admits an open $\Aut(V)$-orbit \cite{MR2822235}. (We will refer to it as the \emph{GLS condition}.)
A variant of a their conjecture, formulated in \cite{MR3866895}, states that for a multisegment $\m$,
$Z(\m)$ is $\square$-irreducible if and only if the GLS condition is satisfied.

Let $S_k$ denote the symmetric group on $k$ elements.
For any permutation $\sigma\in S_k$ consider the multisegment
\[
\m_\sigma=\sum_{i=1}^k[i,2k-\sigma(i)].
\]
A special case of the main result of \cite{MR3866895} is that the representation $Z(\m_\sigma)$ is
$\square$-irreducible if and only if $\sigma$ is smooth, namely $\sigma$ is $3412$ and $4231$ avoiding.
The terminology is motivated by a result of Lakshmibai--Sandhya which asserts that $\sigma$ is smooth
if and only if the corresponding Schubert variety of type $A_{k-1}$ is smooth \cite{MR1051089}.
Another equivalent condition is that the Kazhdan--Lusztig polynomial $P_{e,\sigma}(q)$ with respect to $S_k$ is $1$
\cite{MR788771}.
It turns out that this is also equivalent to the GLS condition in the case at hand, although
we are unaware of a direct relation between the smoothness condition
for Schubert varieties and the GLS condition.

Suppose that $\m_1$ and $\m_2$ are two multisegments.
Using the Arakawa--Suzuki functor \cite{MR1652134}, it is possible to compute the decomposition of
$Z(\m_1)\times Z(\m_2)$ in the Grothendieck group.
The computation involves Kazhdan--Lusztig polynomials with respect to $S_k$ where $k$ is the number
of segments in $\m_1+\m_2$. We refer to  \cite{MR3866895}*{\S10} and the references therein for more details.

In December 2017 I computed all Kazhdan--Lusztig polynomials for each pair of permutations in
$S_{12}$. (In reality, ``only'' about $46\times 10^9$ pairs need to be considered.)
This computation was carried out on a 1 terabyte RAM machine at the Weizmann Institute and resulted
in about $4.3\times 10^9$ different polynomials. More details about this month-long computation (on a single core)
can be found in \cite{1705.06517}.

Following up on this computation, we calculated the decomposition of $\pi\times\pi$ for
any $\pi=Z(\m)$ where $\m$ consists of at most 6 segments.

The most striking feature of the calculation, which was discovered following a question by David Kazhdan,
is that there appears to be a correlation between the length $l$ of $Z(\pi)\times Z(\pi)$ and the minimal
codimension $d$ of an orbit in $\mathcal{C}_{\m}$. We will call the latter the \emph{deficiency} of $\m$.
In particular, the GLS condition is that $d=0$ and by the above, it is conjectured to be equivalent to $l=1$.
(This was verified in all examples checked.)
It is also the case, in all the checked examples that $d=1\iff l=2$ and moreover this condition
is equivalent to $\pi\times\pi$ being multiplicity free.
The only other possible values of $d$ in these examples are $2, 3, 4$ in which cases
the ranges for $l$ are 7--9, 39--53 and 251--257 respectively.

For simplicity, let us say that $\sigma\in S_k$ is \emph{almost smooth} if
$Z(\m_\sigma)\times Z(\m_\sigma)$ has length two.
We say that $\sigma$ is \emph{wild} if it is neither smooth nor almost smooth.

Among the 354 non-smooth permutations in $S_6$, 230 are almost smooth.
For 224 of these, $P_{e,\sigma}(1)=2$, the 6 exceptions being
\[
\sigma=463152, 465132,526413,546213,632541,653421
\]
where $P_{e,\sigma}(1)=3$.

In Table \ref{table1} we list for each almost smooth permutations $\sigma$ in $S_6$
the $w$ such that
\[
Z(\m_\sigma)\times Z(\m_\sigma)=Z(\m_\sigma+\m_\sigma)\oplus Z(\n_w).
\]
Here $w$ is a two-to-one function from $\{1,\dots,12\}$ to $\{1,\dots,6\}$ and
\[
\n_w=\sum_{i=1}^{12}[\lceil\frac i2\rceil,12-w(i)].
\]
In Table \ref{table2} we list for each wild $\sigma\in S_6$, the length $l$ of $\Pi=Z(\m_\sigma)\times Z(\m_\sigma)$,
the number $t$ of distinct constituents in the Jordan--H\"older decomposition of $\Pi$
and the deficiency $d$ of $\m$.

In order to avoid repetition, in both tables we only list one representative for each equivalence class of the equivalence
relation $\sigma\sim\sigma^{-1}\sim\sigma^{w_0}\sim(\sigma^{-1})^{w_0}$
where $\sigma^{w_0}$ is the conjugation of $\sigma$ by the longest permutation $w_0\in S_k$.
(The representation $Z(\m_{\sigma^{-1}})$ is essentially the contragredient of $Z(\m_\sigma)$.
The relation between the decompositions of $Z(\m_\sigma)\times Z(\m_\sigma)$ and
$Z(\m_{\sigma^{w_0}})\times Z(\m_{\sigma^{w_0}})$ is more mysterious but follows from the abovementioned
computation using Kazhdan--Lusztig polynomials.)

In all the cases where $d<4$, i.e., when $\sigma\ne(456123), (562341),(634512)$, any constituent of $\Pi$
occurs with multiplicity one or two.
For $\sigma=(456123)$ where $l=251$, there are 59 constituents occurring with multiplicity one,
81 with multiplicity two, 4 with multiplicity three and 2 with multiplicity 9.
For $\sigma=(562341)$ or $(634512)$ where the length is 257, there are 61 constituents occurring with multiplicity one,
83 with multiplicity two, 4 with multiplicity three and 2 with multiplicity 9.

One can recover from the decomposition of $Z(\m_\sigma)\times Z(\m_\sigma)$, $\sigma\in S_k$,
the decomposition of $Z(\m)\times Z(\m)$ for any multisegment $\m$ consisting of up to $k$ segments.
We will not describe the results obtained in this generality (for $k=6$) but remark that
although more values of $l$ occur in these cases,
the correlation between $l$ and the deficiency of $\m$ continues to hold.


\begin{table}
\begin{tabular}{c|c||c|c||c|c}
$\sigma$ & $w$        & $\sigma$ & $w$        & $\sigma$ & $w$        \\
\hline
125634 & 112235463546 & 126453 & 112246354635 & 135624 & 113325462546 \\
136452 & 113346254625 & 145236 & 112435243566 & 145263 & 112435246635 \\
146253 & 112436254635 & 146352 & 114436253625 & 146523 & 112466352435 \\
153426 & 113524352466 & 153462 & 113524356624 & 153624 & 113524362546 \\
153642 & 113524663524 & 154623 & 115524362436 & 156423 & 112546352436 \\
163542 & 113625463524 & 164352 & 114635243625 & 164532 & 114635463522 \\
215634 & 221135463546 & 216453 & 221146354635 & 235614 & 223315461546 \\
236451 & 223346154615 & 245163 & 221435146635 & 246135 & 221436143655 \\
246153 & 221436154635 & 246351 & 224436153615 & 246513 & 221466351435 \\
251634 & 223511463546 & 253461 & 223514356614 & 253614 & 223514361546 \\
253641 & 223514663514 & 254613 & 225514361436 & 256143 & 221546154633 \\
256314 & 223546351146 & 256413 & 221546351436 & 261453 & 224611354635 \\
263415 & 223614361455 & 263514 & 223615361544 & 263541 & 223615463514 \\
264153 & 224635114635 & 264351 & 224635143615 & 264531 & 224635463511 \\
265341 & 226635143514 & 325614 & 332215461546 & 326451 & 332246154615 \\
346251 & 334426152615 & 351642 & 132514663524 & 352641 & 332514662514 \\
356214 & 332546251146 & 361542 & 132615463524 & 362514 & 332615261544 \\
362541 & 332615462514 & 364251 & 334625142615 & 364521 & 334625462511 \\
365142 & 136625143524 & 365241 & 336625142514 & 365412 & 136655241324 \\
426153 & 241326154635 & 426513 & 241366251435 & 426531 & 241366552413 \\
436512 & 441366251325 & 463152 & 443625113625 & 463251 & 443625132615 \\
463521 & 443625362511 & 465132 & 146635143522 & 465213 & 246635241135 \\
465312 & 146635241325 & 564312 & 154635241326 & 625431 & 261546352413 \\
632541 & 362513462514 & 635421 & 362546352411 & 645321 & 463546352211 \\
653421 & 663524352411 \\
\end{tabular}
\caption{Decomposition of $Z(\m_\sigma)\times Z(\m_\sigma)$ in the almost smooth case}
\label{table1}
\end{table}

\begin{table}
\begin{tabular}{c|c|c|c||c|c|c|c||c|c|c|c}
$\sigma$ & $l$ & $t$ & $d$ & $\sigma$ & $l$ & $t$ & $d$ & $\sigma$ & $l$ & $t$ & $d$ \\
\hline
145623 &  8 & 7 & 2 & 156342 &  8 & 7 & 2 & 163452 &  9 & 8 & 2 \\
245613 &  8 & 7 & 2 & 256134 &  8 & 7 & 2 & 256341 &  8 & 7 & 2 \\
263451 &  9 & 8 & 2 & 264513 &  8 & 7 & 2 & 345612 &  36 & 26 & 3 \\
346152 &  7 & 6 & 2 & 346512 &  8 & 7 & 2 & 351624 &  9 & 8 & 2 \\
354612 &  9 & 8 & 2 & 356124 &  39 & 28 & 3 & 356142 &  7 & 6 & 2 \\
356241 &  8 & 7 & 2 & 356412 &  7 & 6 & 2 & 361452 &  8 & 7 & 2 \\
362451 &  9 & 8 & 2 & 364152 &  9 & 8 & 2 & 364512 &  44 & 31 & 3 \\
426351 &  7 & 6 & 2 & 456123 &  251 & 146 & 4 & 456132 &  9 & 8 & 2 \\
456231 &  44 & 31 & 3 & 456312 &  7 & 6 & 2 & 462351 &  46 & 33 & 3 \\
462513 &  7 & 6 & 2 & 462531 &  7 & 6 & 2 & 463512 &  9 & 8 & 2 \\
465231 &  9 & 8 & 2 & 562341 &  257 & 150 & 4 & 562431 &  7 & 6 & 2 \\
563412 &  53 & 37 & 3 & 563421 &  9 & 8 & 2 & 564231 &  7 & 6 & 2 \\
623451 &  48 & 35 & 3 & 623541 &  7 & 6 & 2 & 624351 &  7 & 6 & 2 \\
624531 &  8 & 7 & 2 & 634521 &  9 & 8 & 2 & 635241 &  9 & 8 & 2 \\
645231 &  53 & 37 & 3 \\
\end{tabular}

\caption{Length, number of distinct constituents of $Z(\m_\sigma)\times Z(\m_\sigma)$
and deficiency in the wild case}
\label{table2}
\end{table}

\clearpage
\newpage

\begin{bibdiv}
\begin{biblist}

\bib{MR1652134}{article}{
      author={Arakawa, Tomoyuki},
      author={Suzuki, Takeshi},
       title={Duality between {$\germ s\germ l_n({\bf C})$} and the degenerate
  affine {H}ecke algebra},
        date={1998},
        ISSN={0021-8693},
     journal={J. Algebra},
      volume={209},
      number={1},
       pages={288\ndash 304},
         url={http://dx.doi.org/10.1006/jabr.1998.7530},
      review={\MR{1652134}},
}

\bib{MR0425031}{article}{
      author={Bern{\v{s}}te{\u\i}n, I.~N.},
      author={Zelevinski{\u\i}, A.~V.},
       title={Induced representations of the group {$GL(n)$} over a {$p$}-adic
  field},
        date={1976},
        ISSN={0374-1990},
     journal={Funkcional. Anal. i Prilo\v zen.},
      volume={10},
      number={3},
       pages={74\ndash 75},
      review={\MR{0425031 (54 \#12989)}},
}

\bib{MR0579172}{article}{
      author={Bernstein, I.~N.},
      author={Zelevinsky, A.~V.},
       title={Induced representations of reductive {${\germ p}$}-adic groups.
  {I}},
        date={1977},
        ISSN={0012-9593},
     journal={Ann. Sci. \'Ecole Norm. Sup. (4)},
      volume={10},
      number={4},
       pages={441\ndash 472},
      review={\MR{0579172 (58 \#28310)}},
}

\bib{MR788771}{article}{
      author={Deodhar, Vinay~V.},
       title={Local {P}oincar\'e duality and nonsingularity of {S}chubert
  varieties},
        date={1985},
        ISSN={0092-7872},
     journal={Comm. Algebra},
      volume={13},
      number={6},
       pages={1379\ndash 1388},
         url={http://dx.doi.org/10.1080/00927878508823227},
      review={\MR{788771 (86i:14015)}},
}

\bib{MR2822235}{article}{
      author={Gei{\ss}, Christof},
      author={Leclerc, Bernard},
      author={Schr{\"o}er, Jan},
       title={Kac-{M}oody groups and cluster algebras},
        date={2011},
        ISSN={0001-8708},
     journal={Adv. Math.},
      volume={228},
      number={1},
       pages={329\ndash 433},
         url={http://dx.doi.org/10.1016/j.aim.2011.05.011},
      review={\MR{2822235}},
}

\bib{MR3314831}{article}{
      author={Kang, Seok-Jin},
      author={Kashiwara, Masaki},
      author={Kim, Myungho},
      author={Oh, Se-jin},
       title={Simplicity of heads and socles of tensor products},
        date={2015},
        ISSN={0010-437X},
     journal={Compos. Math.},
      volume={151},
      number={2},
       pages={377\ndash 396},
         url={http://dx.doi.org/10.1112/S0010437X14007799},
      review={\MR{3314831}},
}

\bib{MR3758148}{article}{
      author={Kang, Seok-Jin},
      author={Kashiwara, Masaki},
      author={Kim, Myungho},
      author={Oh, Se-jin},
       title={Monoidal categorification of cluster algebras},
        date={2018},
        ISSN={0894-0347},
     journal={J. Amer. Math. Soc.},
      volume={31},
      number={2},
       pages={349\ndash 426},
         url={https://doi.org/10.1090/jams/895},
      review={\MR{3758148}},
}

\bib{MR1051089}{article}{
      author={Lakshmibai, V.},
      author={Sandhya, B.},
       title={Criterion for smoothness of {S}chubert varieties in {${\rm
  Sl}(n)/B$}},
        date={1990},
        ISSN={0253-4142},
     journal={Proc. Indian Acad. Sci. Math. Sci.},
      volume={100},
      number={1},
       pages={45\ndash 52},
         url={http://dx.doi.org/10.1007/BF02881113},
      review={\MR{1051089}},
}

\bib{1705.06517}{incollection}{
      author={Lapid, Erez},
       title={Conjectures about certain parabolic {K}azhdan--{L}usztig
  polynomials},
        date={2018},
   booktitle={Geometric aspects of the trace formula},
      series={Simons Symp.},
   publisher={Springer, [Cham]},
       pages={267\ndash 298},
         url={https://doi.org/10.1007/978-3-319-94833-1_9},
        note={arXiv:1705.06517},
}

\bib{MR3866895}{article}{
      author={Lapid, Erez},
      author={M\'{i}nguez, Alberto},
       title={Geometric conditions for {$\square$}-irreducibility of certain
  representations of the general linear group over a non-archimedean local
  field},
        date={2018},
        ISSN={0001-8708},
     journal={Adv. Math.},
      volume={339},
       pages={113\ndash 190},
         url={https://doi.org/10.1016/j.aim.2018.09.027},
      review={\MR{3866895}},
}

\bib{MR0390138}{article}{
      author={Pjasecki{\u\i}, V.~S.},
       title={Linear {L}ie groups that act with a finite number of orbits},
        date={1975},
        ISSN={0374-1990},
     journal={Funkcional. Anal. i Prilo\v zen.},
      volume={9},
      number={4},
       pages={85\ndash 86},
      review={\MR{0390138 (52 \#10964)}},
}

\bib{MR617466}{article}{
      author={Zelevinski{\u\i}, A.~V.},
       title={The {$p$}-adic analogue of the {K}azhdan-{L}usztig conjecture},
        date={1981},
        ISSN={0374-1990},
     journal={Funktsional. Anal. i Prilozhen.},
      volume={15},
      number={2},
       pages={9\ndash 21, 96},
      review={\MR{617466 (84g:22039)}},
}

\bib{MR584084}{article}{
      author={Zelevinsky, A.~V.},
       title={Induced representations of reductive {${\germ p}$}-adic groups.
  {II}. {O}n irreducible representations of {${\rm GL}(n)$}},
        date={1980},
        ISSN={0012-9593},
     journal={Ann. Sci. \'Ecole Norm. Sup. (4)},
      volume={13},
      number={2},
       pages={165\ndash 210},
         url={http://www.numdam.org/item?id=ASENS_1980_4_13_2_165_0},
      review={\MR{584084 (83g:22012)}},
}

\end{biblist}
\end{bibdiv}

\def\cprime{$'$} 

\end{document}